\documentclass{amsart}
\usepackage{amscd,amssymb}

\newcommand{\w}{\omega}

\newcommand{\diam}{\mathrm{diam}\,}
\newcommand{\dist}{\mathrm{dist}}
\newcommand{\mesh}{\mathrm{mesh}\,}
\newcommand{\id}{\mathrm{id}}
\newcommand{\asdim}{\mathrm{asdim}}

\newcommand{\IN}{\mathbb{N}}
\newcommand{\IR}{\mathbb{R}}
\newcommand{\IZ}{\mathbb{Z}}
\newcommand{\IQ}{\mathbb{Q}}
\newcommand{\U}{\mathcal{U}}

\newcommand{\e}{\varepsilon}
\newcommand{\Ra}{\Rightarrow}

\newcommand{\A}{\mathcal A}
\newcommand{\F}{\mathcal F}
\newcommand{\N}{\mathcal N}
\newcommand{\uD}{u\mathcal D}

\newtheorem{mainth}{Theorem}
\newtheorem{maincor}{Corollary}
\newtheorem{theorem}{Theorem}[section]
\newtheorem{lemma}[theorem]{Lemma}
\newtheorem{proposition}[theorem]{Proposition}
\newtheorem{corollary}[theorem]{Corollary}

\newtheorem{mainprob}{Problem}
\newtheorem{problem}[theorem]{Problem}

\begin{document}

\title{The coarse classification of countable abelian groups}
\author{T. ~Banakh}
\address{Instytut Matematyki, Akademia \' Swi\c etokrzyska w Kielcach (Poland),\newline Department of Mathematics, Ivan Franko National University of Lviv (Ukraine)}
\email{tbanakh@yahoo.com}
\author{J.~Higes}
\address{Departamento de Geometr\'{\i}a y Topolog\'{\i}a,
Facultad de CC.Matem\'aticas. Universidad Complutense de Madrid.
Madrid(Spain)}
\email{josemhiges@yahoo.es}
\author{I.~Zarichinyy}
\address{Department of Mathematics, Ivan Franko National University of Lviv (Ukraine)}
\email{ihor.zarichnyj@gmail.com}

%\date{20 August, 2008}

\thanks{ The second named author is supported by Grant AP2004-2494 from the Ministerio de Educaci\' on y Ciencia, Spain and project MEC, MTM2006-0825. 
He thanks Kolya Brodskyi and A. Mitra for helpful discussions. He also thanks Jose Manuel Rodriguez Sanjurjo for his support. Special thanks to Jerzy Dydak 
for all his help and very nice suggestions.}
%\address{}

\begin{abstract} We prove that two countable locally finite-by-abelian groups $G,H$ endowed with proper left-invariant metrics are coarsely equivalent if and only if their asymptotic dimensions coincide and the groups are either both finitely-generated or both are infinitely generated. On the other hand, we show that each countable group $G$ that coarsely embeds into a countable abelian group is locally nilpotent-by-finite. Moreover, the group $G$ is locally abelian-by-finite if and only if $G$ is undistorted in the sense that $G$ can be written as the union $G=\bigcup_{n\in\w}G_n$ of countably many finitely generated subgroups such that each $G_n$ is undistorted in $G_{n+1}$ (which means that the identity inclusion $G_n\to G_{n+1}$ is a quasi-isometric embedding with respect to word metrics on $G_n$ and $G_{n+1}$).
\end{abstract}

\maketitle

%\tableofcontents

\section*{Introduction and main results}

Geometric group theory studies groups viewed as metric spaces, see 
\cite{Gro asym invar}. 
In the case of finitely generated groups the metric to consider is the word metric. 
One of the most important (and still unsolved) problems of this theory is the classification of finitely generated groups up to the quasi-isometry (or, which is the same, to the coarse equivalence), see \cite{Grom},  \cite{Harpe}, \cite{FM2}. %This problem has been solved 
%in some particular cases: for some metabelian groups \cite{Taback} and \cite{FM1}, lamplighter groups \cite{Whyte}, for finitely presented groups of asymptotic dimension $1$, see \cite{Fuji} and \cite{Gentimis}. The general problem is still open even for nilpotent groups.\par

For finitely generated abelian groups this problem has very easy solution:  
a finitely-generated group $G$, being isomorphic to the direct sum of cyclic groups, is coarsely equivalent to the free abelian group $\IZ^{r_0(G)}$ where
$$r_0(G)=\sup\{n\in\w:\mbox{$\IZ^n$ is isomorphic to a subgroup of $G$}\}$$is the {\em torsion free rank} of $G$. So, two finitely generated abelian groups are coarsely equivalent if and only if their torsion free ranks coincide.

We recall that two metric spaces $X,Y$ are ({\em bijectively}) {\em coarsely equivalent} if there are two (bijective) bornologous maps $f:X\to Y$ and $g:Y\to X$ and a real constant $K$ such that $\dist(g\circ f(x),x)\le K$ and $\dist(f\circ g(y),y)\le K$ for all points $x\in X$, $y\in Y$. 
A map $f:X\to Y$ between metric spaces is called {\em bornologous} if for every $\delta<\infty$ the continuity modulus
$$\w_f(\delta)=\sup\{\diam(f(A)):A\subset X,\;\diam(A)\le\delta\}$$is finite.

To treat countable groups as metric spaces we endow them with proper left-invariant metrics. We recall that a metric is {\em proper} if the closed balls with respect to this metric all are compact. By a recent result of Smith \cite{Smith} each countable group $G$ admits a proper left-invariant metric and such a metric is unique up to the coarse equivalence (in the sense that for any two left-invariant proper metrics $d,\rho$ on $G$ the identity map $(G,d)\to (G,\rho)$ is a coarse equivalence). 

Proper left invariant metrics appear naturally when a word metric on a finitely generated group $G$ is restricted to a subgroup $H$ of $G$ (which needs not be  finitely generated).
By studying proper left invariant metrics 
one can extend some of the ideas of classical geometric group theory to countable groups, see \cite{Brod-Dydak-Higes-Mitra}, \cite{BDM}, \cite{Dran-Smith}, \cite{Higes2}, \cite{Sauer}, \cite{Shalom}, \cite{Smith}.  \par

One of important coarse invariants of metric spaces is the asymptotic dimension introduced by M.Gromov. 
Given a metric space $X$ we write $\asdim(X)\le n$ if for every $D\in\IR_+$ there is a cover $\U$ of $X$ such that $$\mesh\U=\sup_{U\in\U}\diam(U)<\infty$$ and $\U$ can be written as the union $\U=\U_0\cup\dots\cup\U_n$ of $(n+1)$ subfamilies $\U_0,\dots,\U_n$ which are {\em $D$-discrete} in the sense that $$\dist(U,V)=\inf\{\dist(u,v):u\in U,\;v\in V\}\ge D$$ for any distinct sets $U,V\in\U_i$, $i\le n$. By definition, the asymptotic dimension $\asdim(X)$ is equal to the smallest integer $n$ such that $\asdim(X)\le n$. If no such an $n$ exists, then $\asdim(X)=\infty$. It well-known that the asymptotic dimension is preserved by coarse equivalences \cite[Ch.9]{Roe}. The asymptotic dimension $\asdim(G)$ of an abelian group $G$ (endowed with a proper left-invariant metric) is equal to the torsion-free rank of $G$, see \cite{Dran-Smith}.

Another property preserved by coarse equivalences is the large scale connectedness. 
A metric space $X$ is defined to be {\em large scale connected} if it is $\e$-connected for some $\e\in\IR_+$. The latter means that any two points $x,y\in X$ can be linked by an $\e$-chain $x=x_0,x_1,\dots,x_n=y$. It is easy to see that a countable group is large scale connected if and only if it is finitely generated.

The main result of this paper is the following classification theorem (see \cite{Higes} for a predecessor of this result). 

\begin{mainth}\label{t1} For two countable abelian groups $G,H$ endowed with proper left-invariant metrics the following conditions are equivalent:
\begin{enumerate}
\item the metric spaces $G,H$ are coarsely equivalent;
\item $\asdim(G)=\asdim(H)$ and the spaces $G,H$ are either both large scale connected or both are not large scale connected;
\item $r_0(G)=r_0(H)$ and the groups $G,H$ are either both finitely-generated or both are infinitely generated.
\end{enumerate}
\end{mainth}

It is interesting to note that the first two conditions in this theorem have metric nature and remain true for metric spaces that are coarsely equivalent to countable abelian groups. For such spaces Theorem~\ref{t1} implies the following classification.

\begin{maincor}\label{c1} If a metric space $X$ of asymptotic dimension $n=\asdim(X)$ is coarsely equivalent to a countable abelian group, then $X$ is coarsely equivalent to
\begin{itemize}
\item $\IZ^n$ iff $X$ is large scale connected;
\item $\IZ^n\oplus(\IQ/\IZ)$ iff $X$ is not large scale connected.
\end{itemize}
\end{maincor}

In light of this result the following problem arises naturally.

\begin{mainprob}\label{prob1} Detect countable groups that are (bijectively) coarsely equivalent to abelian groups.
\end{mainprob}

We shall show that the class of such groups contains all abelian-by-finite groups and all locally finite-by-abelian groups. Let us recall that a group $G$ is $\mathcal P_1$-by-$\mathcal P_2$ where $\mathcal P_1$, $\mathcal P_2$ are two properties of groups, if $G$ contains a normal subgroup $H$ with property $\mathcal P_1$ whose quotient group $G/H$ has the property $\mathcal P_2$. A group
$G$ has a property $\mathcal P$ {\em locally} if each finitely generated subgroup of $G$ has the property $\mathcal P$. 

\begin{mainth}\label{t2} A countable group $G$ is bijectively coarsely equivalent to an abelian group provided $G$ is abelian-by-finite or locally finite-by-abelian.
\end{mainth}

For locally nilpotent undistorted groups this theorem can be reversed.

We define a group $G$ to be {\em undistorted} if $G$ can be written as the union $G=\bigcup_{n\in\w}G_n$ of an increasing sequence $(G_n)_{n\in\w}$ of finitely-generated subgroups such that each group $G_n$ is undistorted in $G_{n+1}$. The latter means that the inclusion $G_n\to G_{n+1}$ is a quasi-isometric embedding with respect to any word metrics on $G_n$ and $G_{n+1}$. 

In Proposition~\ref{p10} we shall show that the class of undistorted groups contains all locally abelian-by-finite groups and also all locally polycyclic-by-finite groups of finite asymptotic dimension.

\begin{mainth}\label{t3} Assume that a countable group $G$ is coarsely equivalent to a countable abelian group. Then
\begin{enumerate}
\item $G$ is locally nilpotent-by-finite;
\item $G$ is locally abelian-by-finite if and only if $G$ is undistorted; 
\item $G$ is locally finite-by-abelian if and only if $G$ is undistorted and locally finite-by-nilpotent.
\end{enumerate}
\end{mainth}

In light of Theorem~\ref{t2}, Theorem \ref{t1} admits the following self-generalization.

\begin{maincor}\label{c2} Two countable locally finite-by-abelian groups $G,H$ are coarsely equivalent if and only if $\asdim(G)=\asdim(H)$ and the groups $G,H$ are either both finitely-generated or else both are infinitely generated.
\end{maincor}

Theorems~\ref{t1} will be proved in Section~\ref{s8} and Theorem~\ref{t3} follows from Proposition~\ref{p5} and Corollary~\ref{c10} below. Problem~\ref{prob1} of detecting groups that are (bijectively) coarsely equivalent to abelian groups will be discussed in Sections~\ref{s9}--\ref{s13}.

\subsection*{Standard notations and conventions} By $\IN$ and $\w$ we denote the sets of positive and non-negative integer numbers, $\IR_+=[0,\infty)$. 

For a metric space $(X,d)$, a point $x_0\in X$, and a real number $r\in\IR_+$ let
$$B_r(x_0)=\{x\in X:d(x,x_0)\le r\}\mbox{ and } O_r(x_0)=\{x\in X:d(x,x_0)<r\}$$denote the closed and open $r$-balls centered at $x_0$.
The product $X\times Y$ of two metric spaces $X,Y$ is endowed with the max-metric
$$d((x,y),(x',y'))=\max\{d_X(x,x'),d_Y(y,y')\}.$$

The neutral element of a group $G$ is denoted by $1_G$. A {\em norm} on a group $G$ is a function $\|\cdot\|:G\to\IR_+$ such that
\begin{enumerate}
\item $\|x\|=0$ if and only if $x=1_G$,
\item $\|x^{-1}\|=\|x\|$,
\item $\|x+y\|\le \|x\|+\|y\|$,
\end{enumerate}
for any $x,y\in G$. 

For any left-invariant metric $d$ on $G$ the formula $\|x\|=d(x,1_G)$ determined a norm on $G$. The metric $d$ can be recovered from this norm by the formula $d(x,y)=\|x^{-1}y\|$.

Due to Smith \cite{Smith} we know that each countable group carries a proper left-invariant metric and such a metric is unique up to the coarse equivalence. Because of that {\em all the groups considered in this paper are countable} (the unique exception appears in the proof of Proposition~\ref{p7} involving connected nilpotent Lie groups). The countable groups are endowed with proper left-invariant metrics. Finitely-generated groups are endowed with word metrics.

\section{Proper left-invariant metrics on homogeneous spaces}

In this section we extend the mentioned result of Smith \cite{Smith} to spaces of the forms $G/H=\{xH:x\in G\}$ where $H$ is a (not necessarily normal) subgroup of a countable group $G$. The space $G/H$ admits a natural left action of the group $G$: $g\cdot xH\mapsto (gx)H$. 

A metric $d$ on $G/H$ is defined to be {\em $G$-invariant} if $d(gx,gy)=d(x,y)$ for all $x,y\in G/H$ and $g\in G$. It is clear that $G/H$ endowed with a $G$-invariant metric is homogeneous as a metric space. 

A metric space $X$ is defined to be {\em homogeneous} if for any two points $x,y\in X$ there is a bijective isometry $f:X\to X$ such that $f(x)=y$. The Baire Theorem guarantees that in a countable proper homogeneous metric space all balls are finite. This simple observation will allow us to prove the uniqueness of proper $G$-invariant metrics on countable spaces $G/H$.

\begin{lemma}\label{l1} Let $H$ be a subgroup of a countable group $G$. For any two proper $G$-invariant metrics $d,\rho$ on $G/H$ the identity map $(G/H,d)\to(G/H,\rho)$ is a coarse equivalence.
\end{lemma}

\begin{proof} It suffices to check that the identity map $(G/H,d)\to(G/H,\rho)$ is bornologous. Observe that for every $\e\in\IR_+$ the real number $$\delta=\max\{\rho(xH,H):d(xH,H)\le\e\}$$is finite because all the balls in $(G/H,d)$ are finite. Then for any cosets $xH,yH\in G/H$ with $d(xH,yH)\le\e$ we get $d(y^{-1}xH,H)=d(xH,yH)\le\e$ by the $G$-invariant property of the metric $d$. Consequently,
$$\rho(xH,yH)=\rho(y^{-1}xH,H)\le\delta,$$
which completes the proof.
\end{proof}

Next, we study the problem of the existence of a proper $G$-invariant metric on $G/H$.
We shall show that such a metric on $G/H$ exists if and only if the subgroup $H$ in quasi-normal in $G$.

 We define a subgroup $H\subset G$ to be 
\begin{enumerate}
\item {\em quasi-normal} if for any $x\in G$ there is a finite subset $F_x\subset G$ such that $x^{-1}Hx\subset F_x H$;
\item {\em uniformly quasi-normal\/} if there is a finite subset $F\subset G$ such that $x^{-1}Hx\subset FH$ for every $x\in G$.
\end{enumerate} 

The (uniform) quasi-normality is tightly related to the usual normality.

\begin{proposition}\label{p1} A subgroup $H$ of a group $G$ is
\begin{enumerate}
\item quasi-normal provided $H$ contains a subgroup of finite index that is normal in $G$;
\item uniformly quasi-normal provided $H$ has finite index in some normal subgroup of $G$.
\end{enumerate}
\end{proposition}

\begin{proof} 1. Assume that $N\subset H$ is a subgroup of finite index, which is normal in $G$. Then for every $x\in G$ we get $x^{-1}Hx\supset x^{-1}Nx=N$. Let $q:G\to G/N$ be the quotient homomorphism. 
Since $N$ has finite index in $H$, it has finite index also in $x^{-1}Hx$.
Then $q(x^{-1}Hx)$ is finite and hence there is a finite subset $F_x\subset G$ such that $x^{-1}Hx\subset F_xN\subset F_xH$.

2. Next, assume that $H$ has finite index in some normal subgroup $N$ of $G$. Then $N=F\cdot H$ for some finite set $F\subset N$ and hence for every $x\in G$ we get $$x^{-1}Hx\subset x^{-1}Nx=N=FH.$$
\end{proof}

The following lemma characterizes subgroups $H\subset G$ whose quotient spaces admit a proper $G$-invariant metric.

\begin{lemma}\label{l2} A subgroup $H$ of a countable group $G$ is quasi-normal in $G$ if and only if the quotient space $G/H$ admits a proper $G$-invariant metric.
\end{lemma}

\begin{proof} To prove the ``if'' part, assume that $d$ is a proper $G$-invariant metric on $G/H$. Since $G/H$ is countable, the Baire Theorem guarantees that all balls with respect to the metric $d$ are finite. In particular, for every $x\in G$ the  ball $B_r(x^{-1}H)$ of radius $r=d(H,x^{-1}H)$ centered at $x^{-1}H$ is finite. So, we can find a finite subset $F_x\subset G$ such that $B_r(x^{-1}H)=\{yH:y\in F_x\}$. We claim that $x^{-1}Hx\subset F_xH$. Take any element $h\in H$ and consider the isometry $f:G/H\to G/H$, $f:yH\mapsto x^{-1}hx\cdot yH$. Observe that $f(x^{-1}H)=x^{-1}H$ and hence $$d(x^{-1}hxH,x^{-1}H)=d(f(H),f(x^{-1}H))=d(H,x^{-1}H)=r.$$ Consequently,
$x^{-1}hxH\in B_r(x^{-1}H)$ and we get the required inclusion $x^{-1}hx\in F_xH$, which means that $H$ is quasi-normal in $G$.

Next, we prove the ``only if'' part. Assume that the subgroup $H$ is quasi-normal in $G$. Fix any proper left-invariant metric $d$ on $G$. Let us show that for any cosets $xH,yH\in G/H$ the Hausdorff distance
$$d_H(xH,yH)=\inf\{\e\in\IR_+:xH\subset O_\e(yH),\; yH\subset O_\e(xH)\}$$is finite.
Here $O_\e(yH)=\{g\in G:\dist(g,yH)<\e\}$ is the open $\e$-neighborhood around $yH$ in $G$.
Consider the point $z=x^{-1}y$ and using the quasi-normality of $H$ in $G$, find a finite subset $F_z\subset G$ such that $z^{-1}Hz\subset F_zH$. Then $Hz\subset zF_zH$ and after inversion, $z^{-1}H\subset HF_z^{-1}z^{-1}$. Find a finite $\e$ such that the open $\e$-ball $O_\e(1_G)$ centered at the neutral element $1_G$ of $G$ contains the finite set $F_z^{-1}z^{-1}$.
Then $$y^{-1}xH=z^{-1}H\subset HF_z^{-1}z^{-1}\subset H\cdot O_\e(1_G)$$and hence
$xH\subset yH\cdot O_\e(1_G)=O_\e(yH)$. 

By analogy, we can check that $yH\subset O_\delta(xH)$ for some $\delta\in\IR_+$.
This proves that the Hausdorff distance is a well-defined metric $d_H$ on $G/H$. It is easy to check that the left-invariant property of the metric $d$ on $G$ implies the $G$-invariant property of the Hausdorff metric $d_H$ on $G/H$.

Also it is easy to see that the quotient map $q:(G,d)\to (G/H,d_H)$ is non-expanding and for every $\e\in\IR_+$ the image $q(O_\e(1_G))$ of the open $\e$-ball centered at the neutral element of $G$ coincides with the open $\e$-ball of $G/H$ centered at $H$.
This observation implies that the metric $d_H$ on $G/H$ is proper.
\end{proof}

Finally, we give conditions on a subgroup $H\subset G$ guaranteeing that the quotient space $G/H$ admits a proper $G$-invariant ultra-metric. 

We shall say that a subgroup $H$ of a group $G$ has {\em locally finite index} in $G$ if $H$ has finite index in each subgroup of $G$, generated by $H\cup F$ for a finite subset $F\subset G$.

\begin{lemma}\label{l3} Let $H$ be a subgroup of a countable group $G$. The space $G/H$ admits a proper $G$-invariant ultra-metric if and only if the subgroup $H$ has locally finite index in $G$.
\end{lemma}

\begin{proof} Assume that the group $H$ has locally finite index in $G$. If $H=G$, then the space $G/H$ is a singleton and trivially admits a proper $G$-invariant ultra-metric. So, we assume that $H\ne G$. Write $G$ as the union $G=\bigcup_{n\in\w}G_n$ of a sequence of subgroups
$$H=G_0\subset G_1\subset G_2\subset\dots
$$such that $H\ne G_1$ and $H$ has finite index in each subgroup $G_n$.
Define a proper $G$-invariant ultra-metric $d$ on $G/H$ letting 
$$d(xH,yH)=\min\{n\in\w: xG_n=yG_n\}.$$
\smallskip

Now assume conversely that the space $G/H$ admits a proper $G$-invariant ultra-metric $\rho$. To show that $H$ has locally finite index in $G$, it suffices to check that for every $r\in\IR_+$ the set $G_r=\{x\in G:\rho(xH,H)\le r\}$ is a subgroup of $G$. (It is clear that $G=\bigcup_{r=1}^\infty G_r$ and the index of $H$ in $G_r$ equals the cardinality of the closed $r$-ball in $G/H$). Take any two point $x,y\in G_r$. Since the ultra-metric $\rho$ is $G$-invariant, we get $$\rho(H,xyH)\le\max\{\rho(H,xH),\rho(xH,xyH)\}=\max\{\rho(H,xH),\rho(H,yH)\}\le r,$$which means that $xy\in G_r$. Also for every $x\in G_r$ we get
$$\rho(H,x^{-1}H)=\rho(xH,xx^{-1}H)=\rho(xH,H)\le r,$$which means that $x^{-1}\in G_r$. Thus $G_r$ is a subgroup of $G$ and we are done.
\end{proof}

Answering Problem 1606 of  \cite{Sanjurjo}, T.Banakh and I.Zarichnyii \cite{BZ} proved that  any two unbounded proper homogeneous ultra-metric spaces are coarsely equivalent. In particular, any such a space is coarsely equivalent to the torsion infinitely-generated group $\IQ/\IZ$.

Since each space $G/H$ endowed with a $G$-invariant metric is homogeneous, we can combine this result with Lemmas~\ref{l3} and \ref{l1} and obtain the following classification.

 \begin{corollary}\label{c3} For any subgroup $H$ of locally finite index in a countable group $G$ the space $G/H$ endowed with a proper $G$-invariant metric is coarsely equivalent to the singleton or to the group $\IQ/\IZ$.
\end{corollary}

The classification of spaces $G/H$ up to the bijective coarse equivalence is more rich. By \cite{BZ}, each homogeneous countable proper ultra-metric space $X$ is bijectively coarsely equivalent to the direct sum
$$\IZ_f=\oplus_{p\in\Pi}\IZ_p^{f(p)}$$
of cyclic groups $\IZ_p=\IZ/p\IZ$ for a suitable function $f:\Pi\to\w\cup\{\infty\}$ defined on the set $\Pi$ of prime numbers. If $f(p)=\infty$ then by $\IZ_p^{f(p)}=\IZ^\infty_p$ we understand the direct sum of countably many copies of the 
group $\IZ_p$. This result combined with Lemmas~\ref{l1} and \ref{l3} implies 

\begin{corollary}\label{c4} For any subgroup $H$ of locally finite index in a countable group $G$ the space $G/H$ endowed with a proper $G$-invariant metric is bijectively coarsely equivalent to the abelian group $\IZ_f$ for a suitable function $f:\Pi\to\w\cup\{\infty\}$.
\end{corollary}

\section{A selection result}

In this section, given a quasi-normal subgroup $H$ of a countable group $G$, we shall study the coarse properties of the quotient map $q:G\to G/H$, $q:x\mapsto xH$. Lemmas~\ref{l1} and \ref{l2} guarantee that $G/H$ carries a proper $G$-invariant metric and such a metric is unique up to the bijective coarse equivalence.

\begin{lemma}\label{l4} For a quasi-normal subgroup $H$ of a countable group $G$ the quotient map $q:G\to G/H$ is bornologous.
\end{lemma}

\begin{proof} By (the proof of) Lemma~\ref{l2}, for any proper left-invariant metric $d$ on $G$  the Hausdorff distance $d_H$ is a proper $G$-invariant metric on $G/H$. It is clear that the quotient map $q:G\to G/H$ is non-expanding (and thus bornologous) with respect to the metrics $d$ and $d_H$. Lemma~\ref{l1} ensures $q:G\to G/H$ is bornologous for any proper $G$-invariant metrics on $G$ and $G/H$.
\end{proof}

Next, we consider the problem of the existence of a bornologous section $s:G/H\to G$ for the quotient map $q:G\to G/H$. A map $s:Y\to X$ is called a {\em section} for a map $f:X\to Y$ if $f\circ s(y)=y$ for all $y\in Y$.

\begin{theorem}\label{p2} Let $H$ be a subgroup of locally finite index in a countable group $G$. For any subsemigroup $S$ of $G$ with $S\cdot H=G$, the quotient map $q:G\to G/H$ has a bornologous section $s:G/H\to S\subset G$.
\end{theorem}

\begin{proof} Since $G$ is countable and $H$ has locally finite index in $G$, we can write the group $G$ as the union $G=\bigcup_{n\in\w}G_n$ of a sequence
$$H=G_0\subset G_1\subset G_2\subset\dots$$
of subgroups of $G$ such that $H$ has finite index in each group $G_n$. If $H\ne G$, then we shall assume that $G_1\ne G_0$. 
Under such a convention the formula
$$\rho(xH,yH)=\min\,\{n\in\w: xG_n=yG_n\}$$ 
determines a proper $G$-invariant ultra-metric on $G/H$.
Fix also any proper left-invariant metric $d$ on $G$.

Let $S\subset G$ be a subsemigroup of $G$ such that $S\cdot H=G$. Attach the unit to $S$ letting $S^1=S\cup\{1_G\}$.

Let $\alpha_0:G_0/H\to S\cap H$ be any map and for every $n\in\IN$ fix a section $\alpha_n:G_n/G_{n-1}\to S^1\cap G_n$ of the quotient map $\pi_n:G_n\to G_n/G_{n-1}$ such that $\alpha_n(G_{n-1})=1_G$ and $\alpha_n(xG_{n-1})\in S\cap xG_{n-1}$ for every $x\in G_n\setminus G_{n-1}$.

Put $s_0=\alpha_0:G_0/H\to S\cap H$ and for every $n\in\IN$ define a section $s_n:G_n/H\to G_n$ of the quotient map $q|G_n:G_n\to G_n/H$ by the recursive formula:
\begin{equation}\label{eq1}
s_n(xH)=\alpha_n(xG_{n-1})\cdot s_{n-1}(\alpha_n(xG_{n-1})^{-1}\cdot xH)
\end{equation}
for $xH\in G_n/H$.

First we show that $s_n(xH)$ is well-defined for every $xH\in G_n/H$. Since $\alpha_n$ is a section of the quotient map $G_n\to G_n/G_{n-1}$, we get $a=\alpha_n(xG_{n-1})\in xG_{n-1}$ and consequently, $x^{-1}a\in G_{n-1}$. Since $H\subset G_{n-1}$, we get $a^{-1}xH\subset a^{-1}xG_{n-1}=a^{-1}x(x^{-1}a)G_{n-1}=G_{n-1}$. Hence $s_{n-1}(a^{-1}xH)$ and $s_n(xH)=a\cdot s_{n-1}(a^{-1}xH)$ are defined.

Observe that for every $x\subset G_{n-1}$, we get $\alpha_n(xG_{n-1})=\alpha_n(G_{n-1})=1_G$ and thus $s_n(xH)=s_{n-1}(xH)$.
This means that $s_n|G_{n-1}/H=s_{n-1}$ and hence we can define a map $s:G/H\to G$ letting $s(xH)=s_n(xH)$ for any $n\in\w$ such that $xH\subset G_n$.

Let us show that the so-defined map  $s$ is a section of $q:G\to G/H$ with $s(G/H)\subset S$.
It suffices to check that for every $n\in\w$ and $x\in G_n$ we get $s(xH)\in S\cap xH$.  This will be done by induction on $n$. If $n=0$, then $s(xH)=s(H)=s_0(H)\in S\cap H$ by the choice of $s_0$. Assume that $s_{n-1}(xH)\in S\cap xH$ for all $x\in G_{n-1}$.
Given any point $x\in G_n\setminus G_{n-1}$ and taking into account that $a=\alpha_n(xG_{n-1})\in S\cap xG_{n-1}$, we get
$$s(xH)=s_n(xH)=a\cdot s_{n-1}(a^{-1}xH)\in a\cdot a^{-1}xH=xH$$and
$$s(xH)=a\cdot s_{n-1}(a^{-1}xH)\subset S\cdot S\subset S$$by the inductive assumption.

Finally, we show that the map $s:G/H\to G$ is bornologous.
Given any $\e\in\IR_+$ find an integer $k\ge \e$ and let 
$$\delta=\diam s(G_k/H).$$
The bornologous property of $s$ will follow as soon as we check the inequality
$d(s(xH),s(yH))\le\delta$ for any $n\in\w$ and points $x,y\in G_n$ with $\rho(xH,yH)\le\e$. This will be done by induction on $n$. If $n\le k$, then $d(s(xH),s(yH))\le\diam s(G_k/H)=\delta$. Assume that the inequality is proved for any $x,y\in G_{n-1}$ with $n>k$. Take two points $x,y\in G_n$ and note that $\rho(xH,yH)\le\e\le k$ implies $xG_k=yG_k$ and hence $xG_{n-1}=yG_{n-1}$. Let $a=\alpha_n(xG_{n-1})=\alpha_n(yG_{n-1})\in xG_{n-1}=yG_{n-1}$ and observe that
$$
\begin{aligned}
d(s(xH),s(yH))&=d(s_n(xH),s_n(yH))=d(a\cdot s_{n-1}(a^{-1}xH),a\cdot s_{n-1}(a^{-1}yH))=\\&=d(s_{n-1}(a^{-1}xH),s_{n-1}(a^{-1}yH))\le\delta.
\end{aligned}$$The last inequality follows from the inductive assumption because $\rho(a^{-1}xH,a^{-1}yH)=\rho(xH,yH)\le\e$.
\end{proof}

\section{Quasi-centralizers and FC-groups}

We recall that the centralizer of a subset $A$ of a group $G$ is the subgroup
$$C(A)=\big\{x\in G:x^A=\{x\}\big\}\mbox{ \ where \ }x^A=\{a^{-1}xa:a\in A\}.$$
By analogy, we define the {\em quasi-centralizer}
$$Q(A)=\{x\in G:\mbox{$x^A$ is finite}\}$$
 of $A$ in $G$.

\begin{lemma} The quasi-centralizer $Q(A)$ of any subset $A\subset G$ is a subgroup of $G$. 
\end{lemma}

\begin{proof} Since $a^{-1}xya=a^{-1}xaa^{-1}ya\subset x^A\cdot y^A$ and thus $(xy)^A\subset x^A\cdot y^A$, we see that $(xy)^A\in Q(A)$ for any $x,y\in Q(A)$. 
Observing that $(x^{-1})^A=(x^A)^{-1}$ we also see that $x^{-1}\in Q(A)$ for each $x\in Q(A)$, which implies that $Q(A)$ indeed is a subgroup of $G$.
\end{proof}

If $A$ is a subgroup of $G$ then we can say a bit more about the subgroup $Q(A)$.
The following lemma easily follows from the definition of the quasi-centralizer.

\begin{lemma} If $H$ is a subgroup of a group $G$, then $x^{-1}Q(H)x=Q(H)$ for every $x\in H$. Consequently, $Q(H)\cdot H=H\cdot Q(H)$ is a subgroup of $G$ containing $Q(H)$ as a normal subgroup.
\end{lemma}

Quasi-centralizers can be used to characterize (locally) FC-groups. 
Following \cite{Baer} we define a group $G$ to be an {\em $FC$-group} if the conjugacy class $x^G$ of each point $x\in G$ is finite.

The subsequent characterization of FC-groups follows immediately from the definitions.

\begin{proposition} A group $G$ is an FC-group if and only if $Q(G)=G$.
\end{proposition}

By an old result of B.H.~Neumann \cite{Neumann} a finitely-generated group is finite-by-abelian if and only if it is an FC-group. This characterization implies the following characterization of locally finite-by-abelian groups.

\begin{proposition} For a group $G$ the following conditions are equivalent:
\begin{enumerate}
\item $G$ is locally FC-group;
\item $G$ is locally finite-by-abelian;
\item $Q(H)\supset H$ for each finitely-generated subgroup $H\subset G$.
\item $Q(H)=G$ for each finitely-generated subgroup $H\subset G$.
\end{enumerate}
\end{proposition}

\section{Bornologous properties of the group operations}

In this section, we establish some bornologous properties of the group operations. 

\begin{lemma}\label{l7} For subsets $A,B$ of a countable group $G$ the multiplication map $$\cdot:A\times B\to G,\;\; \cdot:(a,b)\mapsto ab,$$ is bornologous if and only if $A^{-1}A\subset Q(B)$.
\end{lemma}

\begin{proof} Let $d$ be a  proper left-invariant metric and $\|\cdot\|$ be the corresponding norm on $G$. 

To prove the ``only if'' part, assume that the map $\cdot:A\times B\to G$ is bornologous. Then for every $x\in A^{-1}A$  there is $\delta\in\IR_+$ such that $d(ab,a'b')\le\delta$ for any points $(a,b),(a',b')\in A\times B$ with $\max\{d(a,a'),d(b,b')\}\le \|x\|$. 

Write the point $x\in A^{-1}A$ as $x=a^{-1}_1a_2$ and observe that $d(a_1,a_2)=\|a_1^{-1}a_2\|=\|x\|$. Consequently, for every $b\in B$ we get $$\|b^{-1}xb\|=d(xb,b)=d(a_1^{-1}a_2b,b)=d(a_2b,a_1b)\le\delta.$$
this means that $b^{-1}xb\in B_\delta(1_G)$ and hence $x^B\subset B_\delta(1_G)$, witnessing that $x\in Q(B)$.
\smallskip

Now, assuming that for every $x\in A^{-1}A$ the set $x^B$ is finite, we check that the map $\cdot:A\times B\to G$ is bornologous. Fix any $\e\ge 0$ and consider the finite set $$F=\bigcup\{x^B:x\in A^{-1}A,\; \|x\|\le \e\}.$$ Let $\delta=\max\{\|y\|:y\in F\}$. Now take any two pairs $(a_1,b_1),(a_2,b_2)\in A\times B$ with $\max\{d(a_1,a_2),d(b_1,b_2)\}\le\e$. We claim that $d(a_1b_1,a_2b_2)\le\delta+\e$. Let $x=a_2^{-1}a_1$ and observe that $\|x\|\le\e$ and hence $x^B\subset F$. This yields $\|b_1^{-1}xb_1\|\le\delta$ and
$$\begin{aligned}
d(a_1b_1,a_2b_2)&\le d(a_1b_1,a_2b_1)+d(a_2b_1,a_2b_2)=\\
&=d(b_1^{-1}a_2^{-1}a_1b_1,1)+d(b_1,b_2)\le \|b_1^{-1}xb_1\|+\e\le\delta+\e.
\end{aligned}$$ 
\end{proof}

\begin{lemma}\label{l8} Let $A$ be a subset of a countable group. The inversion operation $$(\cdot)^{-1}:A\to A^{-1},\;\; (\cdot)^{-1}:a\mapsto a^{-1},$$ is bornologous provided  $A^{-1}A\subset Q(A^{-1})$.
\end{lemma}

\begin{proof} Let $d$ be a proper left-invariant metric on the group $G$ and $\|\cdot\|$ be the norm induced by $d$. Given an $\e\in\IR_+$, consider the set $$F=\{a^{-1}xa:a\in A^{-1},\; x\in A^{-1}A,\;\|x\|\le\e\}.$$It is finite because $A^{-1}A\subset Q(A^{-1})$. So $\delta=\max_{x\in F}\|x\|$ is finite. 

The bornologity of the inverse map $(\cdot)^{-1}:A\to A^{-1}$ will follow as soon as we check that $d(x^{-1},y^{-1})\le\delta$ for any points $x,y\in A$ with $d(x,y)\le \e$. It follows that $\|x^{-1}y\|=d(x,y)\le\e$ and hence $yx^{-1}=y(x^{-1}y)y^{-1}\in F$.
Consequently, $d(x^{-1},y^{-1})=\|yx^{-1}\|\le\delta$.
\end{proof}

\section{Two factorization theorems}

In this section, we search for conditions on a countable group $G$ and a quasi-normal subgroup $H\subset G$ guaranteeing that $G$ is bijectively coarsely equivalent to the product $H\times (G/H)$. Here we endow $G/H$ with a proper $G$-invariant metric. Lemmas~\ref{l1} and \ref{l2} guarantee that such a metric on $G/H$ exists and is unique up to the coarse equivalence.

\begin{theorem}\label{t4} Let $H$ be a quasi-normal subgroup of a countable group $G$.  The group $G$ is bijectively coarsely equivalent to $H\times (G/H)$ provided the quotient map $q:G\to G/H$ has a bornologous section $s:G/H\to G$ such that $s(G/H)\subset Q(H)$.
\end{theorem}

\begin{proof} Let $d$ be a left-invariant metric on the group $G$ and $\|\cdot\|$ be the corresponding norm on $G$. Let $s:G/H\to G$ be a bornologous section of the quotient map $q:G\to G/H$ such that $s(G/H)\subset Q(H)$. 

Define a bijective coarse equivalence $f:H\times (G/H)\to G$ by the formula
$f(x,yH)=s(yH)\cdot x$. It is easy to see that the map $f$ is bijective. Taking into account that $Q(H)$ is a group, we conclude that $A^{-1}A\subset Q(H)$ where $A=s(G/H)$.
Now Lemma~\ref{l7} implies that $f$ is bornologous.

It remains to check the bornologous property of the inverse map $$f^{-1}:G\to H\times G/H,\;\;f^{-1}(z)=(s(zH)^{-1}z,q(z)).$$ The bornologity of the quotient map $q$ has been established in Lemma~\ref{l4}. So, it remains to check the bornologity of map $$h:G\to H,\;h:z\mapsto s(zH)^{-1}z.$$ Fix any positive real number $\e$. Since the map $s\circ q:G\to G$ is bornologous, there is $\e_1\in\IR_+$ such that $d(s(zH),s(z'H))\le\e_1$ for any points $z,z'\in G$ with $d(z,z')\le\e$. Since the set $F=\{x^H:x\in Q(H),\;\|x\|\le\e_1\}$ is finite, the number $\delta=\max\{\|y\|:y\in F\}$ is finite too.

Now the bornologous property of the map $h$ will follow as soon as we check that
$d(h(z),h(z'))\le\delta+\e$ for any two points $z,z'\in G$ with $d(z,z')\le \e$. Since $z=s(zH)\cdot h(z)$ and $z'=s(z'H)\cdot h(z')$, we get
$$d(z,z')=d(s(zH)\cdot h(z),s(z'H)\cdot h(z'))=d(s(z'H)^{-1}s(zH)h(z),h(z')).$$
Since $s(G/H)\subset Q(H)$ and $\|s(z'H)^{-1}s(zH)\|=d\big(s(z'H),s(zH)\big)\le\e_1$, we get $h(z)^{-1}s(z'H)^{-1}s(zH)h(z)\in F$ and hence $\|h(z)^{-1}s(z'H)^{-1}s(zH)h(z)\|\le\delta$.
Consequently,
$$
\begin{aligned}
d(h(z),h(z'))&\le d(h(z), s(z'H)^{-1}s(zH)h(z))+d(s(z'H)^{-1}s(zH)h(z),h(z'))=\\
&=\|h(z)^{-1}s(z'H)^{-1}s(zH)h(z)\|+d(z,z')\le \delta+\e.
\end{aligned}
$$
\end{proof}

We recall that a subgroup $H\subset G$ is uniformly quasi-normal if there is a finite subset $F\subset G$ such that $x^{-1}Hx\subset F\cdot H$ for all $x\in G$.

\begin{theorem}\label{t5} Let $H$ be a uniformly quasi-normal subgroup of a countable group $G$ and $A$ be a subset of $G$ such that $A=A^{-1}$ and $Q(A)=G$.  The group $G$ is bijectively coarsely equivalent to $H\times (G/H)$ provided the quotient map $q:G\to G/H$ has a bornologous section $s:G/H\to A\subset G$.
\end{theorem}

\begin{proof} Fix a proper left-invariant metric $d$ on $G$ and a proper $G$-invariant metric $\rho$ on $G/H$. Let $s:G/H\to A$ be a bornologous section of the quotient map $q:G\to G/H$. We claim that the map $$f:H\times G/H\to G,\;\ f:(x,y)\mapsto x\cdot s(y)^{-1},$$
is a bijective coarse equivalence.

Since $Q(A^{-1})=G$, the multiplication and inversion  maps $\mu:H\times A^{-1}\to G$ and $i:A\to A^{-1}$ are bornologous according to Lemmas~\ref{l7} and \ref{l8}.
Then the map $f$ is bornologous as composition of three bornologous maps:
$$H\times G/H\overset{\id\times s}{\longrightarrow} H\times A\overset{\id\times i}\longrightarrow H\times A^{-1}\overset{\mu}\longrightarrow G.$$

It remains to check the bornologity of the inverse map $$f^{-1}:G\to H\times G/H,\;\; f^{-1}:z\mapsto (z\cdot s(z^{-1}H),z^{-1}H).$$
First we check the bornologity of the map $g:G\to G/H$, $g:z\mapsto z^{-1}H$. 

Take any $\e\in\IR_+$. The bornologity of the map $s\circ q:G\to G$ yields a positive real number $\e_1$ such that $d(s(zH),s(z'H))\le\e_1$ for all points $z,z'\in G$ with $d(z,z')\le\e$. It follows from $s(G/H)\subset A=A^{-1}$ and $Q(A)=G$ that the set $F_1=\{axa^{-1}:a\in s(G/H),\;x\in G,\;\|x\|\le\e_1\}$ is finite.
Since the subgroup $H$ is uniformly quasi-normal, there is a finite subset $F_2\subset G$ such that $x^{-1}Hx\subset F_2H$ for all $x\in G$. Finally, let $\delta=\max\{\rho(xH,H):x\in F_1\cdot F_2\}$.

We claim that $\rho(x^{-1}H,y^{-1}H)\le\delta$ for any points $x,y\in G$ with $d(x,y)\le\e$. Observe that $x=s(xH)\cdot h_x$ and $y=s(yH)\cdot h_y$ where $h_x=s(xH)^{-1}\cdot x\in H$ and $h_y=s(yH)^{-1}\cdot y\in H$.
The choice of the number $\e_1$ guarantees that $\|s(yH)^{-1}s(xH)\|=d(s(xH),s(yH))\le\e_1$.

Observe that $$
\begin{aligned}
xy^{-1}&=s(xH)h_xh_y^{-1}s(yH)^{-1}\in s(xH)Hs(yH)^{-1}=\\&=s(yH)\big(s(yH)^{-1}s(xH)\big)s(yH)^{-1}s(yH)Hs(yH)^{-1}\subset\\
&\subset F_1s(yH)Hs(yH)^{-1}\subset F_1F_2H
\end{aligned}$$
and thus $xy^{-1}H=zH$ for some $z\in F_1F_2$. Now the choice of $\delta$ guarantees that
$\rho(x^{-1}H,y^{-1}H)=\rho(H,xy^{-1}H)=\rho(H,zH)\le\delta$.

This completes the proof of the bornologity of the map $g:G\to G/H$, $g:z\mapsto z^{-1}H$ that coincides with the second component of $f^{-1}$. Finally, we check the bornologity of the map $h:G\to H$, $h:z\mapsto z\cdot (s\circ g(z))$, that coincides with the first component of $f^{-1}$. By Lemma~\ref{l7},  the multiplication $\cdot:G\times s(G/H)\to G$, $\cdot:(x,y)\mapsto x\cdot y$, is bornologous. This fact combined with the bornologity of the maps $s$ and $g$ imply the bornologity of the map $h$.
\end{proof}

Now we derive some corollaries from the Factorization Theorems~\ref{t4} and \ref{t5}.

\begin{corollary}\label{c5} Let $H$ be a subgroup of locally finite index in a countable group $G$.
If the subgroup $Q(H)\cdot H$ has finite index in $G$, then $G$ is bijectively coarsely equivalent to $H\times\IZ_f$ for some function $f:\Pi\to\w\cup\{\infty\}$.
\end{corollary}

\begin{proof} First we show that the subgroup $Q=Q(H)\cdot H$ is bijectively coarsely equivalent to $H\times Q/H$. Since $H$ has locally finite index in $Q$, the quotient map $q:Q\to Q/H$ has a bornologous section $s:Q/H\to Q(H)\subset Q$ according to Theorem~\ref{p2}. Applying Theorem~\ref{t4}, we conclude that $Q$ is bijectively coarsely equivalent to $H\times Q/H$. 

The subgroup $Q=Q(H)\cdot H$ has finite index in $G$ and hence is uniformly quasi-normal by Proposition~\ref{p1}. Let $s_1:G/Q\to G$ be any section of the quotient map $q:G\to G/Q$. Taking into account that $G/Q$ is finite, we conclude that $G=Q(s(G/Q)^{-1}\cup s(G/Q))$ and hence $G$ is bijectively coarsely equivalent to $Q\times (G/Q)$ according to Theorem~\ref{t5}. Consequently, $G$ is bijectively coarsely equivalent to $H\times Q/H\times G/Q$. 

By Corollary~\ref{c4},  the space $Q/H$ is bijectively coarsely equivalent to $\IZ_f$ for a suitable function $f:\Pi\to\w\cup\{\infty\}$. The space $G/Q$, being finite, is bijectively coarsely equivalent to the group $\IZ_g$ for a suitable function $g:\Pi\to\w\cup\{\infty\}$. Then the product $Q/H\times G/Q$ is bijectively coarsely equivalent to the product $\IZ_f\times \IZ_g$, which is isomorphic to $\IZ_{f+g}$. Consequently, $G$ is bijectively coarsely equivalent to $H\times \IZ_{f+g}$. 
\end{proof}

%\begin{corollary} A group $G$ belongs to the class $\A_b$ provided $G$ contains a subgroup $H\in\A_b$ of locally finite index in $G$ such that the subgroup $Q(H)\cdot H$ has finite index in $G$.
%\end{corollary}

This corollary implies the ``abelian-by-finite'' part of Theorem~\ref{t2}.

\begin{corollary}\label{c6} Each countable abelian-by-finite group is bijectively coarsely equivalent to an abelian group.
\end{corollary}

\begin{proof} Given a countable abelian-by-finite group $G$, find an abelian subgroup $H$ of finite index in $G$. Then $Q(H)\cdot H\supset H$ has finite index in $G$ and we can apply Corollary~\ref{c5} to conclude that $G$ is bijectively coarsely equivalent to the abelian group $H\times\IZ_f$ for a suitable function $f:\Pi\to\w\cup\{\infty\}$.
\end{proof}

Applying Theorem~\ref{t4} or \ref{t5} to the subgroup $n\IZ$ of the group $\IZ$ we get

\begin{corollary}\label{corZ} For every $n\in\IN$ the group $\IZ$ is bijectively coarsely equivalent to $\IZ\times\IZ_n$.
\end{corollary}

\section{Locally finite-by-abelian groups}\label{s7}

In this section we shall prove the ``locally finite-by-abelian'' part of Theorem~\ref{t2}. This will be done with help of Corollary~\ref{c4} and Selection Theorem~\ref{p2}.\par
The following lemma is rather known in Theory of Groups. We include a proof here for completeness. It is based on the proof of Theorem 1.1 of \cite{Hilton}.

\begin{lemma}\label{laux} The center of a finitely generated finite-by-abelian group has finite index. 
\end{lemma}
\begin{proof}
Let $G = \langle x_1,...,x_n\rangle$ be a finitely generated finite-by-abelian group and let $H$ be a finite subgroup such that $G/H$ is abelian.This means that the commutator $[G, G]$ of $G$ is included in $H$ and it is finite. Now there is a one-one correspondence $y^{-1}\cdot x \cdot y \to x^{-1}\cdot y^{-1}\cdot x \cdot y$  between the set of conjugates of $x$ and some subset of the commutator. Hence $[G:C(x)]$ is finite for all $x \in G$ where $C(x)$ is the centralizer of $x$ in $G$. Therefore $[G: \bigcap_{i=1}^n C(x_i)]$ is finite. It is obvious that $\bigcap_{i=1}^n C(x_i) = Z(G)$, with $Z(G)$ the center of $G$. Thus the center has finite index as claimed.
\end{proof}

\begin{lemma}\label{l9} Each locally finite-by-abelian group $G$ contains a free abelian subgroup $H$ of locally finite index in $G$ such that $Q(H)=G$.
\end{lemma}

\begin{proof} Write $G$ as the union of an increasing sequence $(G_n)_{n\in\w}$ of finitely-generated subgroups of $G$ such that $G_0=\{1_G\}$. Each group $G_n$ is finite-by-abelian, so by Lemma \ref{laux} its center $Z(G_n)$ has finite index in $G_n$. Let $A_0=G_0$ and $A_{n}=A_{n-1}\cdot Z(G_n)$ for $n>0$. It follows that each $A_n$ is an abelian subgroup of finite index in $G_n$. We call a subset $L=\{a_1,\dots,a_n\}$ of an abelian group $A$ {\em linearly independent} if the homomorphism $$h:\IZ^n\to A,\;h:(k_1,\dots,k_n)\mapsto\sum_{i=1}^n k_ia_i,$$is injective.
Let $L_0=\emptyset$ and by induction in each abelian group $A_n$ choose a maximal linearly independent set $L_n$ so that $L_n\supset L_{n-1}$ and $L_n\setminus L_{n-1}\subset Z(G_n)$. Such a choice is possible because for each element $x\in A_n$ there is $k\in\IN$ with $x^k\in Z(G_n)$.

It follows that the subgroup $H_n$  of $A_n$ generated by $L_n$ is isomorphic to $\IZ^{|L_n|}$ while the subgroup $H$ of $G$ generated by the set $L=\bigcup_{n\in\w}L_n$ is a free abelian group. It follows from the choice of the sets $L_i\setminus L_{i-1}\subset Z(G_i)$ that for every $n\in\w$ the set $L\setminus L_{n-1}$ lies in the centralizer $C(G_n)=\{x\in G:\forall y\in G_n\; xy=yx\}$ of the subgroup $G_n$ in $G$. Since $H\cap G_n\supset H_n\cap G_n$ has finite index in $G_n$, the subgroup $H$ has locally finite index in $G$. 

It remains to check that $Q(H)=G$. Take any element $x\in G$ and find $n\in\w$ such that $x\in G_n$. The group $G_n$, being finitely-generated and finite-by-abelian, is an FC-group by \cite{Neumann}. Consequently, the set $x^{G_n}$ is finite and so is the subset $x^{H_n}$ of $x^{G_n}$. We claim that $x^H=x^{H_n}$. Indeed, take any element $h\in H$ and write it as $h=ab$ where $b\in H_n$ and $a$ belongs to the subgroup of $G$ generated by the set $L\setminus L_n$. Since $L\setminus L_n\subset C(G_{n+1})\subset C(G_n)$, the element $a$ commutes with $x$. Consequently, $x^h=h^{-1}xh=b^{-1}a^{-1}xab=b^{-1}xb\in x^{H_n}$.
\end{proof}

Now we are able to prove the ``locally finite-by-abelian'' part of Theorem~\ref{t3}. 

\begin{theorem}\label{t6} Each locally finite-by-abelian group $G$ is bijectively coarsely equivalent to the abelian group $\IZ^m\times\IZ_f$ for some $m\in\w\cup\infty$ and some function $f:\Pi\to\w\cup\{\w\}$.
\end{theorem}

\begin{proof} By Lemma~\ref{l9}, the group $G$ contains a free abelian subgroup $H$ of locally finite index in $G$ such that $Q(H)=G$. The group $H$, being free abelian, is isomorphic to $\IZ^m$ for some $m\in\w\cup\{\infty\}$. By Corollary~\ref{c5}, the group $G$ is bijectively coarsely equivalent to $H\times \IZ_f$ for a suitable function $f:\Pi\to\w\cup\{\infty\}$.
\end{proof}

For groups of infinite asymptotic dimension Theorem~\ref{t6} can be improved as follows.

\begin{theorem}\label{t6a} A countable group $G$ is bijectively coarsely equivalent to $\IZ^\infty$ if and only if $G$ is bijectively coarsely equivalent to an abelian group and $\asdim(G)=\infty$.
\end{theorem}

\begin{proof} The ``only if'' part of this theorem is trivial. To prove the ``if'' part, assume that $\asdim(G)=\infty$ and $G$ is bijectively coarsely equivalent to an abelian group $A$. By Theorem~\ref{t6}, the group $A$ is bijectively coarsely equivalent to $\IZ^m\times\IZ_f$ for some $m\in\w\cup\{\infty\}$ and some function $f:\Pi\to\w\cup\{\infty\}$.
By Corollary 3.3 of \cite{Dran-Smith},
$$\infty=\asdim(G)=\asdim(\IZ^m\times\IZ_f)=\asdim(\IZ^m)+\asdim(\IZ_f)=m+0=m.$$Consequently, $G$ is bijectively coarsely equivalent to 
$$\IZ^\infty\times \IZ_f=\bigoplus_{i=1}^\infty \IZ\times\IZ_{p_i}$$for a suitable sequence $\{p_i\}_{i=1}^\infty\in\{1\}\cup\Pi$. 

It follows from Corollary~\ref{corZ} that for every $i\in \omega$ there is a bijective coarse equivalence $\varphi_i:\IZ\times\IZ_{p_i}\to \IZ$ 
sending the neutral element of $\IZ\times\IZ_{p_i}$ to the neutral element of $\IZ$.
The bijections $(\varphi_i)$ induce a bijective coarse equivalence 
$$\varphi:\oplus_{i=1}^\infty \IZ\times\IZ_{p_i}\to\oplus_{i=1}^\infty \IZ=\IZ^\infty,\;\;\varphi(g_i)\mapsto(\varphi_i(g_i)).$$

Consequently, 
$$G\sim A\sim\IZ^\infty\times\IZ_f\simeq\oplus_{i=1}^\infty\IZ\times\IZ_{p_i}\sim \IZ^\infty$$
where $\simeq $ means isomorphic and $\sim$ bijectively coarsely equivalent.
\end{proof}

\section{The proof of Classification Theorem~\ref{t1}}\label{s8}

Given two countable abelian groups $G,H$ endowed with proper left-invariant metrics we need to prove the equivalence of the following three conditions:
\begin{enumerate}
\item the metric spaces $G,H$ are coarsely equivalent;
\item $\asdim(G)=\asdim(H)$ are the spaces $G,H$ are either both large scale connected or both are not large scale connected;
\item $r_0(G)=r_0(H)$ and the groups $G,H$ are either both finitely-generated or both are infinitely-generated.
\end{enumerate}

The implication $(1)\Ra(2)$ follows immediately from the invariantness of the asymptotic dimension and the large scale connectedness under coarse equivalences.

The equivalence $(2)\Leftrightarrow(3)$ follows from the subsequent two lemmas. The first of them was proved in \cite{Dran-Smith}.

\begin{lemma}[Dranishnikov, Smith]\label{l10} The asymptotic dimension $\asdim(G)$ of a countable abelian group $G$ equals the torsion-free rank of $G$.
\end{lemma}

The second lemma also should be known. We give a proof for completeness.

\begin{lemma}\label{l11} A countable group $G$ endowed with a proper left-invariant metric is large scale connected if and only if $G$ is finitely generated.
\end{lemma}

\begin{proof}
If $G$ is a finitely generated group then its word metric is large-scale connected, in fact it is $\e$-connected for $\e=1$. For the converse suppose that a countable group $G$ endowed with a proper left invariant metric $d$ is $\e$-connected for some $\e\in\IR_+$. As $d$ is proper the closed ball $B_\e(1_G)$ is finite. Since $(G,d)$ is $\e$-connected, every element $g\in G$ can be linked with the unit $1_G$ by a chain $1_G = g_0,g_1,  ..., g_n = g$ such that $d(g_i, g_{i+1}) \le  \e$ for all $i\le n$. By the left invariant condition this implies 
$g_i^{-1}\cdot g_{i+1} \in B_\e(1_G)$ and as we can write $g = (g_0^{-1}\cdot g_1) \cdot (g_1^{-1} \cdot g_2)\cdots(g_{n-1}^{-1}\cdot g_n)$ we get that $G$ is generated 
by the finite set $B_\e(1_G)$. 
\end{proof}

Finally, the implication $(2)\Ra(1)$ of Theorem~\ref{t1} follows from the ``abelian'' version of Corollary~\ref{c1}.

\begin{lemma} A countable abelian group $G$ of asymptotic dimension $n=\asdim(G)$ is coarsely equivalent to
\begin{itemize}
\item $\IZ^n$ if $G$ is finitely-generated;
\item $\IZ^n\times(\IQ/\IZ)$ if $G$ is infinitely generated.
\end{itemize}
\end{lemma}

\begin{proof} By Lemma~\ref{l10}, $r_0(G)=\asdim(G)=n$.

 If $G$ is finitely generated, then $G$ is isomorphic to $\IZ^{n}\oplus T$ for some finite group $T$. Since the projection $\IZ^n\oplus T\to\IZ^n$ is a coarse equivalence, the group $G$ is coarsely equivalent to $\IZ^n$.

If $G$ is infinitely-generated, then by Theorem~\ref{t6}, $G$ is bijectively coarsely equivalent to the group $\IZ^m\oplus \IZ_f$ for some $m\in\w\cup\{\infty\}$ and some function $f:\Pi\to\w\cup\{\infty\}$. By Corollary 3.3 in \cite{Dran-Smith}, $$n=\asdim(G)=\asdim(\IZ^m\oplus\IZ_f)=\asdim(\IZ^m)+\asdim(\IZ_f)=m+0=m.$$

If the group $\IZ_f$ is infinitely-generated, then it is coarsely equivalent to $\IQ/\IZ$ according to Corollary~\ref{c3}. In this case $G$ is bijectively coarsely equivalent to $\IZ^n\times (\IQ/\IZ)$.

If $\IZ_f$ is finitely generated, then it is finite and then $n=\infty$ because otherwise the coarsely equivalent groups $\IZ^n\times\IZ_f$ and $G$ would be finitely generated. Since the abelian group $\IZ^n\times\IZ_f$ has infinite asymptotic dimension we can apply Theorem~\ref{t6a} to conclude that $\IZ^n\times\IZ_f$ is bijectively coarsely equivalent to any other countable abelian group of infinite asymptotic dimension, in particular, to the group $\IZ^n\times(\IQ/\IZ)$.
\end{proof}

\section{Coarse embeddings into abelian groups}\label{s9}

From now on we will be occupied with Problem~\ref{prob1} of detecting groups that are coarsely equivalent to abelian groups. 
We start with studying properties of countable groups that coarsely embed into countable abelian groups. 

We recall that a map $f:X\to Y$ between two metric spaces is called a {\em coarse embedding} if $f:X\to f(X)$ is a coarse equivalence, where $f(X)$ is endowed with the metric induced from $Y$.

\begin{proposition}\label{p5} If a countable group $G$ coarsely embeds into a countable abelian group, then $G$ is locally nilpotent-by-finite.
\end{proposition}

\begin{proof} This proposition is an easy corollary of the famous theorem of M.Gromov \cite{Gromov2}.

\begin{theorem}[Gromov]\label{t7} A finitely generated group $G$ is nilpotent-by-finite if and only if $G$ has polynomial growth.
\end{theorem}

We recall a group $G$ with finite generating set $S=S^{-1}$ has {\em polynomial growth} if there are constants $C$ and $d$ such that for every $n\in\IN$ we get $|S^n|\le C\cdot n^d$. 

To prove Proposition~\ref{p5}, assume that $f:G\to A$ is a coarse embedding of a countable group $G$ into an abelian group $A$. We need to prove that each finitely-generated subgroup of $G$ is nilpotent-by-finite. We lose no generality assuming that $G$ is finitely-generated. Fix a finite generating subset $S=S^{-1}$ for the group $G$ and consider the word metric $d$ corresponding to this set $S$.

It follows that the metric space $(G,d)$ is {\em $\e$-connected} for $\e=1$. Since $f$ is bornologous, the number $$r=\w_f(1)=\sup\{\diam f(B):B\subset G,\;\diam(B)\le 1\}$$ is finite. Let $B$ stands for the closed $r$-ball in the abelian group $A$, centered at the neutral element of $A$. Since the countable group $A$ is endowed with a proper left-invariant metric, the ball $B$ is finite. The subgroup of $A$, generated by the set $B$ is abelian and consequently, has polynomial growth. So we can find constants $C,d$ such that $|B^n|\le Cn^d$ for all $n\in\IN$. 

Since $f$ is a coarse embedding, the number $m=\sup\{|f^{-1}(a)|:a\in A\}$ is finite. Then for every $n\in\IN$ we get
$$|S^n|\le m\cdot|f(S^n)|\le m\cdot |B^n|\le m\cdot C\cdot n^d,$$
which means that the group $G$ has polynomial growth and hence is nilpotent-by-finite by the Gromov's Theorem~\ref{t7}.
\end{proof}

It turns out that for finitely-generated groups Proposition~\ref{p5} can be reversed.

\begin{proposition}\label{p6} A finitely-generated group $G$ coarsely embeds into a countable abelian group if and only if $G$ is nilpotent-by-finite.
\end{proposition}

\begin{proof} The ``only if'' part of this proposition follows from Proposition~\ref{p5}. To prove the ``if'' part, take any finitely-generated nilpotent-by-finite group $G$ and fix a finite generating set $S=S^{-1}$ for $G$. By the Bass Theorem 2 in \cite{Bass},  the group $G$ has polynomial growth of some integer degree $d$ in the sense there is a constants $C$ such that for every $n\in\IN$ we get
$$C^{-1}n^d\le |S^n|\le Cn^d.$$
This fact implies that the word metric $\rho$ on $G$ induced by the generating set $S$ is {\em doubling} in the sense that there is a constant $C\in\IR_+$ such that for every $r\in\IR_+$ and every $x_0\in G$ we get
$$|B_{2r}(x_0)|\le C\cdot |B_r(x_0)|,$$
where $B_r(x_0)$ denotes the closed $r$-ball centered at $x_0$.

It is easy to check that the formula
$$\sqrt{\rho}(x,y)=\sqrt{\rho(x,y)},\;\;x,y\in G,$$ determines a proper left-invariant metric $\sqrt{\rho}$ on $G$. Consequently, the identity map $(G,\rho)\to (G,\sqrt{\rho})$ is a coarse equivalence.

Since the metric $\rho$ on $G$ is doubling, the Assouad Embedding Theorem \cite{As}  yields a bi-Lipschitz embedding $f:(G,\sqrt{\rho})\to \IR^m$. Since the identity map $(G,\rho)\to (G,\sqrt{\rho})$ is a coarse equivalence, the  map $f:G\to \IR^m$ seen as a map from $(G,\rho)$ to $\IR^m$ is a coarse embedding.

It is clear that the map $[\cdot]:\IR\to\IZ$, assigning to each real number $x\in\IR$ its integer part $[x]$, is a coarse equivalence which induces a coarse equivalence
$$i:\IR^m\to\IZ^m,\;\;[\cdot]^m:(x_1,\dots,x_m)\mapsto([x_1],\dots,[x_m]).
$$

Now we see that the composition $g=i\circ f:G\to \IZ^m$ is a required coarse embedding.
\end{proof}

\section{Quasi-isometric embeddings into abelian groups}

In the previous section we characterized finitely-generated groups admitting coarse embeddings into countable abelian groups as nilpotent-by-finite groups. 
In this section we shall obtain a similar characterization for quasi-isometric embeddings.

We recall that a map $f:X\to Y$ is called a {\em quasi-isometric embedding} if there are constants $L,C$ such that for every $x,x'\in X$ we get
$$\frac1Ld_X(x,x')-C\le d_Y(f(x),f(x'))\le Ld_X(x,x')+C.$$
Here the right-hand inequality means that the map $f$ is {\em asymptotically Lipschitz}.

A quasi-isometric embedding $f:X\to Y$ is called a {\em quasi-isometry} if  $f(X)$ is large in $Y$ is the sense that there is a real constant $K$ such that for every $y\in Y$ there is $x\in X$ with $d_Y(y,f(x))\le K$. 

It is easy to see that each asymptotically Lipschitz map is bornologous. Consequently, each quasi-isometry is a coarse equivalence. The converse is true if the spaces are roughly geodesic.

We define a metric space $X$ to be {\em roughly geodesic} if there is a constant $C$ such that any two points $x,y\in X$ can be linked by a chain $x=x_0,\dots,x_n=y$ such that $n\le d(x,y)+1$ and $d_X(x_{i-1},x_i)\le C$ for all $i\le n$.

The following easy lemma is well-known.

\begin{lemma}\label{l13} Each bornologous map $f:X\to Y$ defined on a roughly geodesic metric space is asymptotically Lipschitz.
\end{lemma}

This lemma implies that each coarse equivalence between roughly geodesic metric spaces is a quasi-isometry. In contrast, a coarse embedding between roughly geodesic spaces needs not be a quasi-isometric embedding.  Here is a simple counterexample.

Consider the Heisenberg group $UT_3(\IZ)$ consisting of unitriangular matrices of the form 
$$
H_{(x,y,z)}=\left(\begin{array}{ccc} 
1 & x & y\cr
0 & 1 & z\cr
0 & 0 & 1
\end{array}\right)\mbox{ where }x,y,z\in\IZ.$$
According to \cite[1.38]{Roe}, for any word metric  on $UT_3(\IZ)$ there are positive constants $c,C$ such that 
the distance $\|H_{(x,y,z)}\|$ from a matrix $H_{x,y,z}$ to the identity matrix lies in the interval
$$c(|x|+|y|+\sqrt{|z|})\le\|H_{(x,y,z)}\|\le C(|x|+|y|+\sqrt{|z|}).$$
This means that the homomorphism 
$$h:\IZ\to UT_3(\IZ),\; h:z\mapsto H_{(0,0,z)},$$
is a coarse embedding but not a quasi-isometric embedding.

The Heisenberg group is the simplest example of a non-abelian nilpotent torsion-free group. For such groups we have the following characterization.

\begin{proposition}\label{p7} A finitely-generated nilpotent torsion-free group $G$ is abelian if and only if $G$ admits a quasi-isometric embedding into a finitely-generated abelian group.
\end{proposition}

\begin{proof} This proposition will be derived from its ``continuous'' version proved by Scott Pauls in \cite{Pauls}.

\begin{theorem}[Pauls] If a connected nilpotent Lie group $G$ endowed with a left-invariant Riemannian metric admits a quasi-isometric embedding into a finite-dimensional Euclidean space, then $G$ is abelian.
\end{theorem}

The reduction of Proposition~\ref{p7} to the Pauls' Theorem will be made with help of a classical result of Malcev \cite{Mal}, see also \cite[2.18]{Rag}.

\begin{theorem}[Malcev] Each finitely-generated torsion-free nilpotent group is isomorphic to a uniform lattice in a simply-connected nilpotent Lie group.
\end{theorem}

We recall that a discrete subgroup $H$ of a topological group $G$ is called a {\em uniform lattice} if the quotient space $G/H$ is compact.
\medskip

In order to prove Proposition~\ref{p7}, take any quasi-isometric embedding $f:G\to A$ of a finitely-generated torsion-free nilpotent group $G$ into a finitely-generated abelian group $A$. Since $A$ is quasi-isometric to $\IZ^m$ for $m=r_0(G)$, we lose no generality assuming that $A=\IZ^m$.

By the Malcev Theorem, the group $G$ can be considered as a uniform lattice in a connected nilpotent Lie group $L$. Endow $L$ with a proper left-invariant Riemannian metric on $L$. Since the quotient group $L/G$ is compact, there is a compact subset $K\subset L$ such that $G\cdot K=L$. For each point $x\in L$ we can choose a point $g(x)\in G$ such that $x\in g(x)\cdot K$. It follows that the map $g:L\to G$ is a quasi-isometric embedding and so is the composition $g\circ f:L\to A=\IZ^m\subset\IR^m$.
Applying the Pauls' Theorem, we conclude that the Lie group $L$ is abelian and so is its subgroup $G$.
\end{proof}

\begin{corollary}\label{c7} A finitely-generated group $G$ quasi-isometrically embed into a finitely-generated abelian group if and only if $G$ is abelian-by-finite.
\end{corollary}

\begin{proof} If a finitely-generated group $G$ is abelian-by-finite, then $G$ contains an abelian subgroup $A$ of finite index. In this case the inclusion $A\subset G$ is a quasi-isometry and hence $G$ quasi-isometrically embeds into the finitely-generated abelian group $A$.
This proves the ``if'' part of the corollary.

To prove the ``only if'' part, assume that  a finitely-generated group $G$ quasi-isometrically embeds into a finitely-generated abelian group. Since each quasi-isometric embedding is a coarse embedding, the group $G$ is nilpotent-by-finite according to Proposition~\ref{p5}. Consequently, $G$ contains a nilpotent subgroup $N$ of finite index in $G$. This subgroup is finitely-generated because $G$ is finitely generated. Since $N$ has finite index in $G$ the inclusion map $i:N\to G$ is a quasi-isometry.

Let $T$ be the set of periodic elements of the nilpotent group $N$. By Theorem 17.2.2 of \cite{Kar}, $T$ is a finite normal subgroup of $N$ and the quotient group $N/T$ is a nilpotent finitely-generated torsion-free group. Since $T$ is finite, the quotient homomorphism $q:N\to N/T$ is a quasi-isometry. Taking into account that $q:N\to N/T$ and $i:N\to G$ are quasi-isometries, we conclude that the group $N/T$ is quasi-isometric to $G$ and hence admits a quasi-isometric embedding into a finitely-generated abelian group. By Proposition~\ref{p7}, the group $N/T$ is abelian. Consequently, the group $N$ is finite-by-abelian and hence its center $Z(N)$ has finite index in $N$ by Lemma \ref{laux}. This means that $N$ is abelian-by-finite and so is the group $G$.
\end{proof}

With help of Corollary~\ref{c7} we can characterize finite-by-abelian groups as follows.

\begin{corollary}\label{c8} A finitely-generated group is finite-by-abelian if and only if it is abelian-by-finite and finite-by-nilpotent.
\end{corollary}

\begin{proof} The ``only if'' part follows from Lemma \ref{laux}.

To prove the ``if'' part, assume that a finitely-generated group $G$ is abelian-by-finite and finite-by-nilpotent. Let $A\subset G$ be an abelian subgroup of finite index and $F\subset G$ be a finite normal subgroup with nilpotent quotient $G/F$. Let $q:G\to G/F$ be the quotient homomorphism. 
It follows that the nilpotent group $G/F$ is finitely-generated. Consequently, the torsion subgroup $T$ of $G/F$ is finite. Replacing the subgroup $F$ by $F\cdot q^{-1}(T)$, we can assume that the quotient group $G/F$ is torsion-free. Since the maps $A\to G$ and $G\to G/F$ are quasi-isometries, we see that the nilpotent torsion-free group $G/F$ is quasi-isometric to the abelian group $A$. By Proposition~\ref{p7}, the group $G/F$ is abelian. Consequently, the group $G$ is finite-by-abelian.
\end{proof}

\begin{corollary}\label{c9} A countable group is locally finite-by-abelian if and only if it is locally abelian-by-finite and locally finite-by-nilpotent.
\end{corollary}

\section{Undistorted groups} 

In this section we establish some properties of undistorted groups. We recall that a finitely-generated subgroup $H$ of a finitely generated group $G$ is called {\em undistorted} in $G$ if the inclusion map $i:H\to G$ is a quasi-isometric embedding with respect to some (equivalently, any) word metrics on the groups $H,G$, see \cite{Mitra}. An example of a distorted subgroup is the cyclic subgroup generated by the matrix $H_{(0,0,1)}$ in the Heisenberg group $UT_3(\IZ)$.

We say that a subgroup $H$ of a countable group $G$ is a {\em bornologous retract} of $G$ if there is a bornologous map $r:G\to H$ such that $r(h)=h$ for all $h\in H$. If 
$r$ can be chosen to be a group homomorphism, then we say that $H$ is {\em complemented} in $G$.

\begin{lemma}\label{l14} A finitely-generated subgroup $H$ of a finitely-generated group $G$ is undistorted in $G$ provided one of the following conditions holds:
\begin{enumerate}
\item $H$ is a bornologous retract of $G$;
\item $H$ is complemented in $G$;
\item $H$ has finite index in $G$;
\item $G$ is abelian-by-finite;
\item $G$ is polycyclic-by-finite and $\asdim(H)=\asdim(G)$;
\end{enumerate}
\end{lemma}

\begin{proof} Fix word metrics on the groups $H,G$. Since the word metric on $H$ is roughly geodesic, the identity inclusion $i:H\to G$ is asymptotically Lipschitz by Lemma~\ref{l13}.

1. Assume that $r:G\to H$ is a bornologous retraction. Since the word metric on $G$ is roughly geodesic, $r$ is asymptotically Lipschitz and so is its restriction $r|i(H)=i^{-1}:i(H)\to H$, witnessing that $i$ is a quasi-isometric embedding.

2. If $H$ is complemented in $G$, then $H$ is undistorted in $G$, being a bornologous 
retract of $G$.

3. If $H$ has finite index in $G$, then the inclusion $H\to G$ is a coarse equivalence and hence a quasi-isometry because the word metrics on $H$ and $G$ are roughly geodesic.

4. Assume that $G$ is abelian-by-finite and let $A$ be an abelian subgroup of finite index in $G$. Since $G$ is finitely-generated, so are the abelian subgroups $A$ and $A\cap H$. Since $A$ has finite index in $G$, the subgroup $A\cap H$ has finite index in $H$. Let $L_1$ be a maximal linearly independent subset of the abelian group $A\cap H$. It can be enlarged to a maximal linearly independent subset $L$ in $A$. It follows that the free abelian subgroup $F_1$ generated by $L_1$ is complemented in the free abelian group $F$, generated by the set $L$ in $A$. Consequently, $L_1$ is undistorted in $L$. 

By the maximality of $L_1$ and $L$ the subgroup $F_1$ has finite index in $A\cap H$ while $L$ has finite index in $A$.
Consequently, the identity inclusions $L_1\to A\cap H\to H$ and $L\to A\to G$ are quasi-isometries. Since $L_1$ is undistorted in $L$, the subgroup $H$ undistorted in $G$.

5. Assume that $G$ is polycyclic and $\asdim(H)=\asdim(G)$. By \cite{Dran-Smith}, $\asdim(G)$ equals to the Hirsh rank $hr(G)$ of $G$, while $\asdim(H)$ is equal to $hr(H)$. By induction of the Hirsh length it is easy to check that the equality $hr(H)=hr(G)$ implies that $H$ has finite index in $G$ and hence $H$ is undistorted in $G$ by the third item of the lemma.
\end{proof}

A  group $G$ is defined to be {\em undistorted\/} provided $G$ is the union of a non-decreasing sequence $(G_n)_{n\in\w}$ of finitely-generated subgroups such that each subgroup $G_n$ is undistorted in $G_{n+1}$.

\begin{proposition}\label{p10} A countable group $G$ is undistorted provided  one of the following conditions holds:
\begin{enumerate}
\item $G$ contains a finitely-generated subgroup $H$ of locally finite index in $G$;
\item $G$ is locally abelian-by-finite;
\item $G$ is locally polycyclic-by-finite and $\asdim(G)<\infty$;
\item $G$ coarsely embeds into a countable abelian group and $\asdim(G)<\infty$.
\end{enumerate}
\end{proposition}

\begin{proof} 1,2. The first two items of this proposition follow immediately from the items 3,4 of Lemma~\ref{l14}.

3. Assume that the group $G$ is locally polycyclic-by-finite and $\asdim(G)<\infty$. By \cite{Dran-Smith}, the group $G$ contains a finitely-generated subgroup $H$ such that $\asdim(H)=\asdim(G)$. Write $G$ as the union of a sequence $H=G_0\subset G_1\subset\cdots$ of finitely-generated subgroups. Each group $G_n$ is polycyclic-by-finite because $G$ is locally polycyclic-by-finite. It follows from $\asdim(H)=\asdim(G)$ that $\asdim(G_n)=\asdim(H)$ for all $n$. Now Lemma~\ref{l14}(5) implies that each group $G_n$ is undistorted in $G_{n+1}$, which means that the group $G$ is undistorted.

4. Assume that $G$ coarsely embeds into a countable abelian group and $\asdim(G)<\infty$.
By Proposition~\ref{p5}, the group $G$ is locally nilpotent-by-finite. Since finitely-generated nilpotent groups are polycyclic, $G$ is locally polycyclic-by-finite and hence undistorted according to the preceding item.
\end{proof}

\section{Groups that are coarsely equivalent to abelian groups}

In this section we study undistorted groups that are coarsely equivalent to abelian groups. The following theorem is our main result in this direction.

\begin{theorem}\label{t10} If a countable undistorted group $G$ is coarsely equivalent to a countable abelian group $A$, then $G$ is locally abelian-by-finite.
\end{theorem}

\begin{proof} Since $G$ is coarsely equivalent to $A$ there are two bornologous maps $f:G\to A$, $g:A\to G$ and a constant $r\in\IR_+$ such that $d_G(g\circ f(x),x)\le r$ and $d_A(f\circ g(y),y)\le r$ for all $x\in G$ and $y\in A$.

The group $G$, being undistorted, can be written as the union $G=\bigcup_{n\in\w}G_n$ of a non-decreasing sequence $(G_n)_{n\in\w}$ of finitely-generated subgroups of $G$ such that each subgroup $G_n$ is undistorted in $G_{n+1}$. We shall also assume that the group $G_0$ contains the closed $r$-ball $B_r(1_G)\subset G$.

Since the restrictions $f|G_n:G_n\to A$ are coarse embeddings, the groups $G_n$ are nilpotent-by-finite according to Proposition~\ref{p5}. We claim that every group $G_n$ is abelian-by-finite. Find a finitely-generated abelian subgroup $A_n\subset A$ such that $f(G_n)\subset A_n$. The abelian-by-finite property of the group $G_n$ will follow from Corollary~\ref{c7} as soon as we check that the restriction $f|G_n:G_n\to A_n$ is a quasi-isometric embedding with respect to word metrics on $G_n$ and $A_n$. First note that this map is asymptotically Lipschitz because the word metric on $G_n$ is roughly geodesic. 

Since the subgroup $A_n$ is finitely generated and $g:A\to G$ is bornologous, the image $g(A_n)$ lies in some finitely-generated subgroup of $G_m$ of $G$. Since the word metric on $A_n$ is roughly geodesic, the map $g|A_n:A_n\to G_m$ is asymptotically Lipschitz and so is its restriction $g|f(G_n):f(G_n)\to G_m$. Taking into account the definition of $r$ and the inclusion $B_r\subset G_0\subset G_n$, we conclude that  $g\circ f(G_n)\subset G_n$.
Since the inclusion $G_n\subset G_m$ is a quasi-isometric embedding, the map $g:f(G_n)\to G_n$ is asymptotically Lipschitz with respect to the word metric on the group $G_n$. This implies that $f|G_n:G_n\to A_n$ is a quasi-isometric embedding. So, we can apply Corollary~\ref{c7} to conclude that the group $G_n$ is abelian-by-finite and consequently, the group $G$ is locally abelian-by-finite.
\end{proof}

Unifying Theorem~\ref{t10} with Proposition \ref{p10} and Corollary~\ref{c9} we obtain the last two items of Theorem~\ref{t3} (the first item of this theorem follows from Proposition~\ref{p5}).

\begin{corollary}\label{c10} If a countable group $G$ is coarsely equivalent to an abelian group, then
\begin{enumerate}
\item $G$ is locally abelian-by-finite if and only if $G$ is undistorted;
\item $G$ is finite-by-abelian if and only if $G$ is undistorted and locally finite-by-nilpotent.
\end{enumerate}
\end{corollary}

Another corollary concerns groups of finite asymptotic dimension.

\begin{corollary} If a countable group $G$ of finite asymptotic dimension is coarsely equivalent to an abelian group, then $G$ is locally abelian-by-finite.
\end{corollary}

\begin{proof} By Proposition~\ref{p10}(4), the group $G$ is undistorted and by Theorem~\ref{t1}, $G$ is locally abelian-by-finite.
\end{proof}

Finally, we characterize finitely generated groups that are (bijectively) coarsely equivalent to abelian groups.

\begin{corollary} For a finitely generated group $G$ the following conditions are equivalent:
\begin{enumerate}
\item $G$ is abelian-by-finite.
\item $G$ is coarsely equivalent to a countable abelian group;
\item $G$ is bijectively coarsely equivalent to a countable abelian group;
\end{enumerate}
\end{corollary}

\begin{proof} The implication $(1)\Ra(3)$ was proved in Corollary~\ref{c6} while $(3)\Ra(2)$ is trivial. Since finitely generated groups are undistorted, the implication $(2)\Ra(1)$ follows from Theorem~\ref{t10}.
\end{proof}

\section{Some Open Problems}\label{s13}

In this section we discuss some open problems related to the topic of the paper. First we draw a diagram displaying the interplay between various classes of groups close to being abelian in a suitable (coarse) sense. 

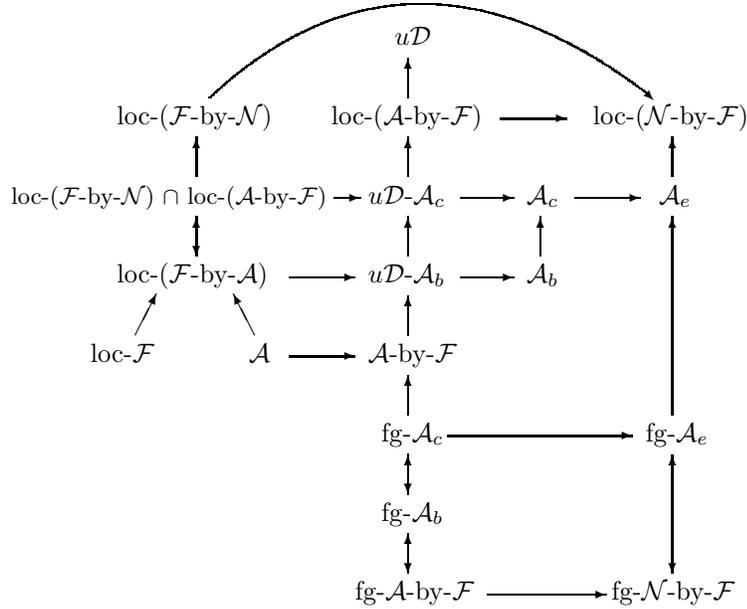
\begin{picture}(300,255)(40,-10)
\put(180,0){fg-$\A$-by-$\F$}
\put(200,15){\vector(0,1){10}}\put(200,20){\vector(0,-1){10}}
\put(230,2){\vector(1,0){45}}

\put(277,0){fg-$\N$-by-$\F$}
\put(300,20){\vector(0,1){35}}\put(300,40){\vector(0,-1){30}}

\put(190,30){fg-$\A_b$}
\put(200,45){\vector(0,1){10}}\put(200,50){\vector(0,-1){10}}

\put(190,60){fg-$\A_c$}
\put(200,70){\vector(0,1){15}}
\put(215,62){\vector(1,0){70}}
\put(290,60){fg-$\A_e$}
\put(300,70){\vector(0,1){75}}

\put(80,90){loc-$\F$}
\put(97,100){\vector(1,2){8}}
\put(140,90){$\A$}
\put(142,100){\vector(-1,2){8}}
\put(155,92){\vector(1,0){25}}
\put(185,90){$\A$-by-$\F$}
\put(200,100){\vector(0,1){15}}

\put(90,120){loc-($\F$-by-$\A$)}
\put(120,135){\vector(0,1){10}}\put(120,140){\vector(0,-1){10}}
\put(155,122){\vector(1,0){25}}

\put(185,120){$\uD$-$\A_b$}
\put(200,130){\vector(0,1){15}}
\put(220,122){\vector(1,0){20}}
\put(245,120){$\A_b$}
\put(250,130){\vector(0,1){15}}

\put(50,150){\small loc-($\F$-by-$\N$) $\cap$ loc-($\A$-by-$\F$)}
\put(120,160){\vector(0,1){15}}
\put(172,152){\vector(1,0){10}}

\put(185,150){$\uD$-$\A_c$}
\put(200,160){\vector(0,1){15}}
\put(220,152){\vector(1,0){20}}
\put(245,150){$\A_c$}
\put(263,152){\vector(1,0){25}}
\put(295,150){$\A_e$}
\put(300,160){\vector(0,1){15}}

\put(90,180){loc-($\F$-by-$\N$)}
\qbezier[0](125,190)(200,260)(290,192)
\put(288,194){\vector(1,-1){5}}
\put(171,180){loc-($\A$-by-$\F$)}
\put(200,190){\vector(0,1){15}}
\put(235,182){\vector(1,0){25}}
\put(270,180){loc-($\N$-by-$\F$)}

\put(195,210){$\uD$}
\end{picture}

In this diagram we consider the following classes of countable groups:
\begin{itemize}
\item $\A$ of abelian groups;
\item $\A_b$ of groups that are bijectively coarsely equivalent to abelian groups;
\item $\A_c$ of groups that coarsely equivalent to abelian groups;
\item $\A_e$ of groups that coarsely embed into abelian groups;
\item $\N$ of nilpotent groups;
\item $\F$ of finite groups;
\item $\mathrm{fg}$ of finitely-generated groups;
\item $\uD$ of undistorted groups.
\end{itemize}
The other classes are obtained as intersections, extensions and localizations of those classes.

In light of this diagram the following questions arise naturally.

\begin{problem}\label{prob2} Are the classes $\A_b$ and $\A_c$ equal?
\end{problem}

\begin{problem}\label{prob3} Is each group $G\in\A_c$ locally abelian-by-finite (equivalently, undistorted)?
\end{problem}

\begin{problem}\label{prob4} Is each locally abelian-by-finite group (bijectively) coarsely equivalent to an abelian group?
\end{problem}

\begin{problem}\label{prob5} Does each locally nilpotent-by-finite group coarsely embed into an abelian group?
\end{problem}

The answers to Problems~\ref{prob2}, \ref{prob4} and \ref{prob5} are not known even for groups of finite asymptotic dimension.

Now let us discuss some possible counterexamples to Problems~\ref{prob3} and \ref{prob4}.
The first group that comes to mind is the union
$$UT_\infty(\IZ)=\bigcup_{n=1}^\infty UT_n(\IZ)$$ of the chain
$$UT_1(\IZ)\subset UT_2(\IZ)\subset \cdots$$ of groups of unitriangular matrices, where each group $UT_{n-1}(\IZ)$ is identified with a subgroup of $UT_{n}(\IZ)$ consisting of unitriangular matrices $(a_{i,j})_{i,j=1}^n$ with $a_{i,n}=0$ for all $i<n$.
The group $UT_\infty(\IZ)$ is locally nilpotent and has infinite asymptotic dimension.

\begin{problem} Does the group $UT_\infty(\IZ)$ admit a coarse embedding into an abelian group? Is $UT_\infty(\IZ)$ (bijectively) coarsely equivalent to an abelian group?
Is $UT_\infty(\IZ)$ undistorted?
\end{problem}

Another interesting concrete group to consider in the wreath product $\IZ\wr S_\infty$ where $S_\infty=\bigcup_{n=1}^\infty S_n$ is the group of finitely supported bijections of $\w$. The group $\IZ\wr S_\infty$ is a semidirect product $\IZ^\infty\ltimes S_\infty$ with respect to the coordinate permutating action of $S_\infty$ on the direct sum $\IZ^\infty$ of countably many infinite cyclic groups. The group $\IZ\wr S_\infty$ is locally abelian-by-finite but not locally finite-by-abelian.

\begin{problem} Is $\IZ\wr S_\infty$ coarsely equivalent to an abelian group? Does $\IZ\wr S_\infty$ coarsely embed into an abelian group?
\end{problem}

Next, we consider some hereditary properties of the classes $\A_b$, $\A_c$, and $\A_e$.
It is clear that those classes are closed under taking direct products. In fact, Theorems~\ref{t4} and \ref{t5} imply that those classes are closed under certain more general product constructions. 

\begin{proposition} A countable group $G$ belongs to the class $\A_i$ for $i\in\{b,c,e\}$ provided $G=A\cdot B$ for two subgroups $A,B\in \A_i$ of $G$ such that $A\cap B=\{1_G\}$ and one of the following conditions holds:
\begin{enumerate}
\item $A$ is quasi-normal in $G$ and $B\subset Q(A)$;
\item $A$ is uniformly quasi-normal in $G$ and $G=Q(B)$.
\end{enumerate}
\end{proposition}

Another stability property of the classes $\A_b$, $\A_c$ and $\A_e$ follows from Corollary~\ref{c5}.

\begin{proposition} A countable group $G$  belongs to the class $\A_i$ for $i\in\{b,c,e\}$ if $G$ contains a subgroup $H\in\A_i$ of locally finite index in $G$ such that the group $Q(H)\cdot H$ has finite index in $G$.
\end{proposition}

It is clear that the class $\A_e$ is closed under taking subgroups.

\begin{problem} Are the classes $\A_b$ and $\A_c$ closed under taking subgroups? 
\end{problem}

\begin{problem} Are the classes $\A_b$, $\A_c$, and $\A_e$ closed under taking quotient groups? 
\end{problem}

We recall that a class of countable groups is {\em local} provided a countable group $G$ belongs to the class if and only if each finitely generated subgroup of $G$ belongs to that class.

\begin{problem}
Are the classes $\A_b$, $\A_c$, and $\A_e$ local? 
\end{problem}

Let us observe that the last question has affirmative solution if and only if Problems~\ref{prob4} and \ref{prob5} have affirmative solutions.

\end{document}